\documentclass[a4paper, 10pt, conference]{ieeeconf}      

\IEEEoverridecommandlockouts                              
\usepackage{graphics} 
\usepackage{epsfig} 
\usepackage{amsmath} 
\usepackage{placeins}
\usepackage[labelformat=simple]{subcaption}

\title{\LARGE \bf
Lagrangian formulation for mixed traffic flow including two-wheelers*
}

\author{Sosina Gashaw$^{1}$, J\'{e}r\^{o}me H\"{a}rri$^{1}$, Paola Goatin$^{2}$
\thanks{$^{1}$Sosina Gashaw and J\'{e}r\^{o}me H\"{a}rri are with EURECOM, 06904 Sophia Antipolis, France
        {\tt\small }}
\thanks{$^{2}$Paola Goatin is with Inria Sophia Antipolis M\'{e}diterran\'{e}e, Universit\'e C\^ote d’Azur, Inria, CNRS, LJAD, 06904 Sophia Antipolis,France
        {\tt\small }}%
}

\begin{document}

\maketitle
\thispagestyle{plain}
\pagestyle{plain}

\begin{abstract}
Lagrangian formulation of kinematic wave provides a more accurate representation than the most commonly used Eulerian formulation. Furthermore, Lagrangian representation offers a flexibility to study certain traffic phenomena (e.g. capacity drop, traffic delay, trajectory). The natural resemblance of Lagrangian representation to road traffic data collection methods (probe vehicles) also renders it suitable for applications such as traffic motoring. This paper presents a multi-class Lagrangian representation for a traffic flow consisting of cars and two-wheelers. A numerical method for solving the continuity equation is presented. We compare the results obtained with Eulerian and Lagrangian formulations. 

\end{abstract}

\section{INTRODUCTION}
The integration of powered two wheelers (PTWs) to intelligent transport systems as well as the development of PTWs specific innovative transport solutions depends upon the understanding of their mobility behaviors and interaction with other road users. However, PTWs create a peculiar traffic flow effects that cannot be reproduced by the currently available models. \par

The flow of vehicles can be modeled at different granularities. Macroscopic models study collective behavior whereas microscopic representation model individual vehicle mobility. Mesoscopic models exhibit both macroscopic and microscopic behaviors, combine the aggregate level modeling in macroscopic approach with specific individual vehicle characteristics such as probabilistic lane changing and turning ratio. The choice of the modeling approach depends on different factors such as the required level of detail, accuracy, efficiency. Macroscopic modeling is an efficient and a preferable approach for studying analytical flow properties, since it allows to establish a closed form relationship between flow variables.  \par
Macroscopic traffic flow models most commonly apply the kinematic wave theory developed by Lighthill, Whitham and Richards \cite{c6,c7} (LWR). In order to integrate traffic heterogeneity (vehicle and driver), LWR model is extended to multi-class flow model. The variation among vehicle class is expressed in relation to maximum speed, perception difference to area/space occupancy, total/effective density. Multi-class LWR models are usually solved in the Eulerian coordinates. In Eulerian formulation, the evolution of flow properties such as density, flow, speed, etc., are evaluated at fixed points. However, recent studies show that Lagrangian representation, which tracks the evolution of flow properties of vehicle/platoon of vehicles, offers several advantages over the Eulerian representation, with the main benefits being numerical accuracy \cite{c2} and flexibility to easily incorporate traffic phenomena (e.g. capacity drop \cite{c15}) and vehicle characteristics.\par

In Lagrangian systems, the LWR model is formulated in $(N,t)$ coordinate system \cite{c2}. Cumulative vehicle count (N) is found to be more suitable for certain traffic flow analysis \cite{c5,c10} and also makes it easier to establish a connection between follow-the-leader and LWR models \cite{c3}. For a mixed traffic of cars and trucks, the Lagrangian formulation is given in \cite{c1}, which identifies vehicle classes with a specific fundamental diagram, jam density, wave speed. However, the interaction between vehicle classes is disregarded. Similarly, in \cite{c4} another Lagrangian representation for multi-class LWR model is proposed. Nonetheless, these models are intended to characterize mixed traffic of cars and trucks. The discretization schemes fall short of describing correctly multi-class flows that have different characteristics from cars and trucks mixed flow, for example mixed flow of cars and two-wheelers, thus requiring further modification. \par

To this end, in this paper, we propose a Lagrangian formulation for traffic flow consisting of cars and two-wheelers (PTWs). The derivation follows the Eulerian multi-class LWR model in \cite{c8}, where the fundamental diagram and the parameters for the fundamental diagram are defined uniquely for each class, and are also adapted to the traffic condition.  We provide a discretization method applicable for solving any type of multi-class LWR model, including cars and PTWs mixed flow. Moreover, we propose an approach to reproduce a follow-the-leader type behavior using the Lagrangian representation. In car traffic, there is an ordered type of flow where the $n^{th}$ vehicle (follower) follows the ${(n-1)}^{th}$ vehicle. PTWs usually do not respect such ordered flow. To accurately represent the abreast movement of two-wheelers in Lagrangian representation, we introduce sub-lanes. We test the equivalence of the Lagrangian and Eulerian representation.

From application standpoint, the Lagrangian representation is convenient to analyze vehicle-specific data such as trajectories and travel times. Using the spacing and speed data collected from probe vehicles together with the traffic flow model formulated in Lagrangian coordinate, traffic state can be estimated accurately \cite{c11}. Moreover, in hybrid traffic flow models, Lagrangian model is used in conjunction with Eulerian representations \cite{c13}. \par

The rest of the paper is organized as follows. First, a formulation for traffic flow consisting PTWs and cars in Eulerian and the Lagrangian approaches is discussed. Thereafter, a discretization technique is presented. Following the Numerical examples and discussion, we wind up by giving concluding remarks. \par
\section{Lagrangian Formulation of Multi-class LWR}
We first introduce the Eulerian representation of the mixed flow of cars and PTWs, and then we show the transformation to Lagrangian coordinates. Multi-class LWR models distinguish the characteristics of each of the vehicle class. Different methods have been applied to represent accurately the distinctive features exhibited, depending on the involved vehicle types. In this study, our interest is in modeling mixed cars and PTWs flow. Thus, a model developed in \cite{c8} is used as a reference model. The model is based on free space distribution, wherein the difference in vehicle size, lateral and longitudinal gap acceptance and maximum speed are the factors that differentiate vehicle classes. 
The continuum equation which holds for each class is written as 
\begin{equation}
 \frac{\partial \rho_i(x,t)}{\partial t} + \frac{\partial q_i(x,t)}{\partial x}=0, \qquad \qquad  i=1,2,
\end{equation}
where $\rho_i$ and $q_i$ denote density and flow of class i, respectively, over space $x$ and time $t$. Class specific flow, speed and density are related by the equation
\begin{equation} \label{eq:MLWRq}
q_i(x,t)=\rho_i(x,t) v_i(x,t),\qquad \qquad i=1,2.
\end{equation}
The speed $v_i$ for the individual vehicle class $i$ is a function of the densities of both classes and derived based on the assumption that the flow of vehicles is dictated by available free space \cite{c8}:
\begin{equation}\label{eq:speed_fn}
v_i=V_i(\rho_1,\rho_2)= v_i^f \left(\int_{r^c_i}^\infty f(l(\rho_1,\rho_2))) \, \mathrm{d}l\right)
\end{equation}
where $v_i^f, r^c_i$ and $f(l(\rho_1,\rho_2)))$ stand for the maximum speed, critical lateral gap and the probability density function of free space distribution, respectively. \par
\begin{figure}[thpb]
\centering
\captionsetup{justification=centering}
\begin{subfigure}{0.75\linewidth}
\includegraphics[width=\textwidth]{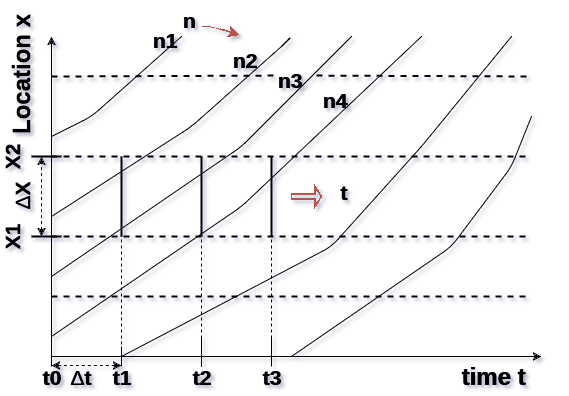}
\caption{Eulerian fixed frame}
\label{fig:euler}
\end{subfigure}
\begin{subfigure}{0.75\linewidth}
\includegraphics[width=\textwidth]{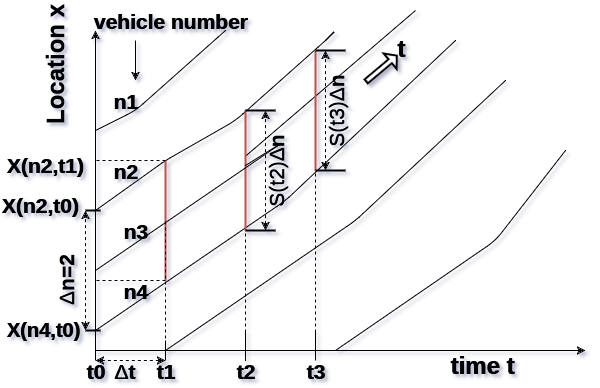}
\caption{Lagrangian moving frame}
\label{fig:lag}
\end{subfigure}

\caption{A schematic of Lagrangian and Eulerian approach}
\end{figure}

 
The Eulerian representation describes the evolution of the traffic state variables at a fixed point in space (Figure \ref{fig:euler}). Whereas, the Lagrangian view deal with the flow properties observed along the trajectory of vehicles (Figure \ref{fig:lag}). \par

The mathematical form of the conservation law in Lagrangian coordinates depends on the chosen coordinate system. Here, we take $(n,t)$ coordinate system. Moreover, there are two methods that are used to represent multi-class flows in Lagrangian coordinates.  In the first method, there are separate Lagrangian coordinates for each vehicles class (method 1). On contrary, in the second method (method 2) there is one Lagrangian reference frame that moves with one of the selected vehicle class. Thus, for the other vehicle classes the conservation equation is derived based on this Lagrangian reference frame. In a situations where tracking of each vehicle is needed (e.g for class specific controls \cite{c9}) method 1 is suitable. Otherwise, method 2 is a computationally efficient approach, for instance to investigate the impact of PTWs on cars flow, or vice versa.  

\subsection{Method 1}
By taking spacing as a state variable, the conservation equation in $(n,t)$ coordinate system is written as \cite{c2}:
\begin{equation}\label{mt1_a}
\frac{\partial s_i(x(t),t)}{\partial t} + \frac{\partial v_i(n,t)}{\partial n}=0 \qquad \qquad i=1,2
\end{equation}
\begin{equation}
s=\frac{-\partial x}{\partial n},~~\rho=\frac{-\partial n}{\partial x}=1/s 
\end{equation}
where $s$ and $v$ denote, respectively, the average spacing and the speed associated to a group of vehicle/s labeled $n$. Vehicle groups are labeled in time order. The conversation equation applies for each vehicle class. Moreover, the grouping of vehicle and the labeling of vehicle groups is done separately for each vehicle class. This representation also take an assumption that vehicles in a group neither disband or merge with other group. Class specific speed-spacing fundamental relation has the following form:
\begin{equation}\label{mt1_b}
v_i=V(s_1,s_2)
\end{equation}

\begin{figure}[thpb]
\captionsetup{justification=centering}
\begin{subfigure}{0.45\linewidth}
\includegraphics[width=\textwidth]{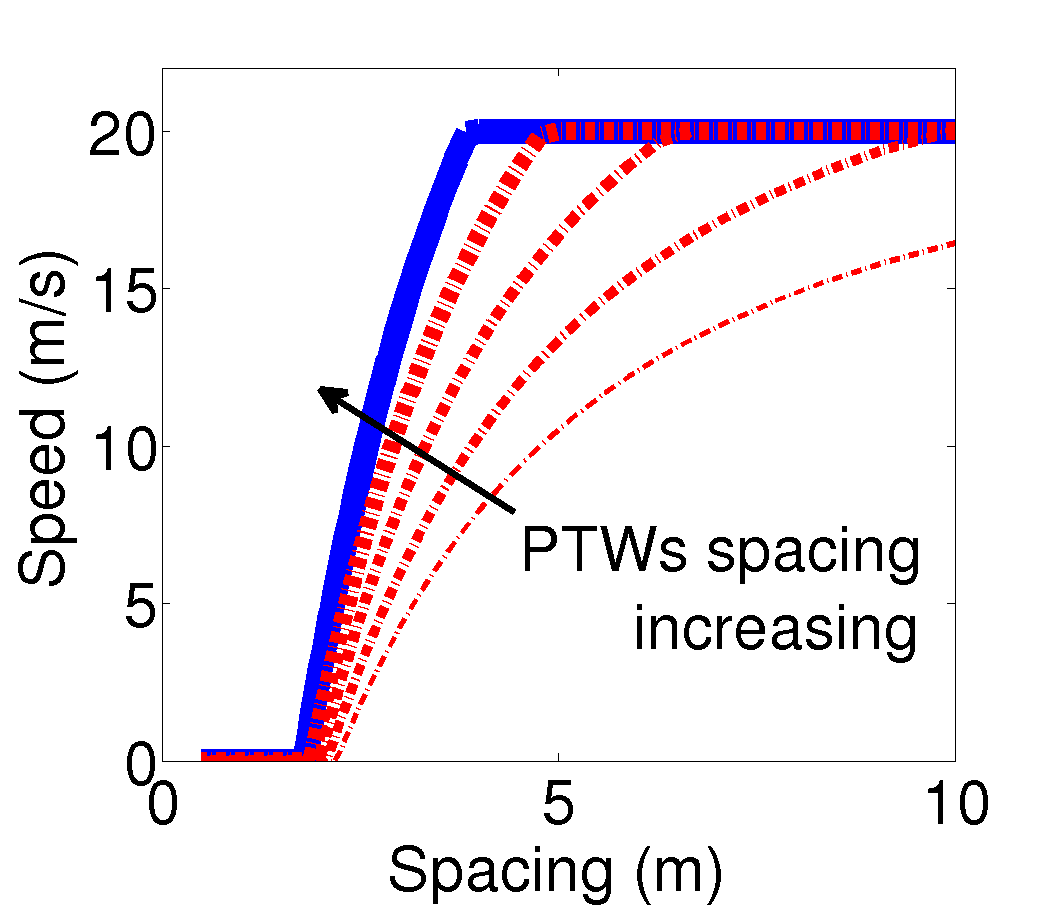}
\caption{}
\label{fig:fd1}
\end{subfigure}
\begin{subfigure}{0.45\linewidth}
\includegraphics[width=\textwidth]{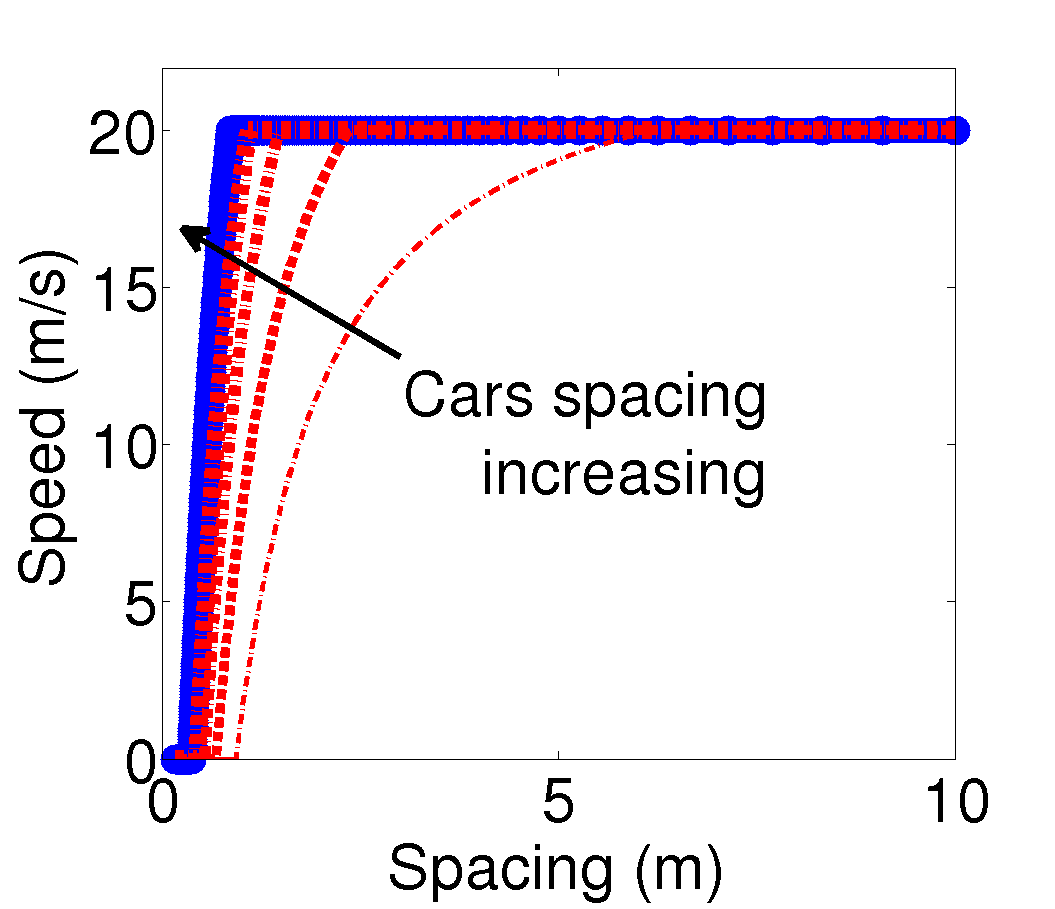}
\caption{}
\label{fig:fd2}
\end{subfigure}
\caption{Speed-spacing fundamental diagram (a) for Cars (b) for PTWs, $V_{1,max}=V_{2,max}=20m/s$ }
\label{fig:fd}
\end{figure}

Speed-space fundamental diagram (FD) for PTWs and cars is given in Figure \ref{fig:fd}. As illustrated in the figure, the fundamental diagram for each class changes with the spacing/density of the other vehicle class. 
\subsection{Method 2}
In the above multi-class Lagrangian conservation equation, individual vehicle class has a separate labeling (cumulative vehicle count). \cite{c4} proposed an alternative formulation, where the Lagrangian coordinates move with a reference vehicle class and only vehicle of this class are counted. In another word, the evolution of traffic state variable of the carrier vehicle class and other vehicle classes being carried inside is tracked. \par
The motion of the reference (carrier) class is governed by:
\begin{equation}\label{mt2_a}
\frac{\partial s_r(x(t),t)}{\partial t} + \frac{\partial v_r(n,t)}{\partial n}=0,
\end{equation}
For the rest vehicle classes:
\begin{equation}\label{mt2_b}
\frac{\partial s_r/s_i}{\partial t} + \frac{\partial \left((v_r-v_i)/s_i\right)}{\partial n}=0,
\end{equation}
or equivalently it can be formulated in non conservative form
\begin{equation*}
\frac{\partial s_i}{\partial t} + \frac{s_i}{s_r} \frac{\partial v_i}{\partial n}- \frac{v_i-v_r}{s_r}\frac{\partial s_i}{\partial n}=0 
\end{equation*}
 where the subscript $r$ and $i$ refers to, respectively, the reference vehicle class and other vehicle classes.\par
In the conservation equation given above, the traffic state variable are spacing (s) and speed (v). When density ($\rho$) is used instead of spacing, the equation takes the following conservation form.
\begin{equation}
\frac{\partial (1/\rho_1)}{\partial t} + \frac{\partial v}{\partial n}=0 
\end{equation}
\begin{equation}
\frac{\partial (\rho_i/\rho_1)}{\partial t}+ \frac{\partial (\rho_i (v_1-v_i))}{\partial n}=0
\end{equation}
where $\rho_1 > 0$ always. 
\section{Discretization scheme}
We apply the following numerical scheme to find the solution of Eq. (\ref{mt1_a}) - (\ref{mt1_b}) (Method 1). The $n$ domain is subdivided into $\Delta n$ sized clusters of vehicles (cells). An approximation of the average spacing $s$ over each cluster is updated at each time step $\Delta t$. Applying Godunov scheme, the numerical solution of the conservation equation is approximated by
\begin{equation} \label{eq:m1}
s_i^{t+\Delta t}=s_i^{t} -\frac{\Delta t}{\Delta n}(V_{i+1/2}-V_{i-1/2})
\end{equation}
where $V_{i+1/2}$ and $V_{i-1/2}$ are the fluxes (speeds) at the boundaries of cell $i$.
\begin{equation*}
V_{i+1/2}=V(s_{1,i},s_{2,i},...),~ V_{i-1/2}=V(s_{1,i-1},s_{2,i-1},...)
\end{equation*}
Therefore, equation \ref{eq:m1} becomes
\begin{equation}
s_i^{t+\Delta t}=s_i^{t} -\frac{\Delta t}{\Delta n}(V(s_{1,i},...)-V(s_{1,i-1},...) 
\end{equation}
which is similar to the direct difference approximation of the conservation equation. To obtain a stable solution $\Delta t$ should be restricted to Courant-Friedrichs-Lewy (CFL) condition, i.e. 
\begin{equation*}
\Delta t \leq \frac{\Delta n}{\max(\lambda)}
\end{equation*}
where $\lambda$ is the wave speed. Following the definition of the flux (speed) at the boundary, the trajectory (location) $X$ \cite{c1} of each cluster can be updated using 
\begin{equation}
X(i,t+\Delta t)= X(i,t) + \Delta t*V(s_{1,i},s_{2,i},...)
\end{equation}
\begin{figure}
\centering
\includegraphics[width=0.8\linewidth]{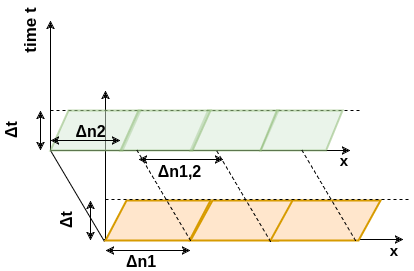}
\caption{n-t domain discretization, separate coordinate for each vehicle class}
\label{fig:vc_disc}
\end{figure}
For each vehicle class, clusters don not overlap each others. However, clusters of different vehicles class may overlap or occupy the same position. For example, in Fig. \ref{fig:vc_disc} the first cluster of vehicle class $1$ overlaps with two clusters of the other vehicle class. To compute $V(s_{1,i},s_{2,i})$, we need to approximate $s_{2,i}$ value in cluster i of vehicle class $1$.
\begin{equation*}
s_{2,i}^{(1)}=\frac{\Delta n_1 s_{1,i}}{\int_{X(i)}^{X(i-1)}\frac{1}{s_2(x)}dx} 
\end{equation*}
where $s_2(x)$ is a function describing the average spacing $s$ of class $2$ as a function of location $x$. For the general case, 
\begin{equation}
s_{c,i}^{(j)}=\frac{\Delta n_j s_{j,i}}{\int_{X(i)}^{X(i-1)}\frac{1}{s_c(x)}dx} ~~ c=1,2,...
\end{equation}
where j and c denote, respectively, the vehicle class cluster $i$ belongs to and the other vehicle classes. For $j=v$, the integration is reduced to $\Delta n_j$, thus $s_{v,i}^{(j)}=\Delta s_{j,i}$.\par

The discretization method introduced above is for method 1, where each of the vehicle class is counted and grouped separately. However, there is an alternative representation (method 2) as shown in Eq. (7)-(8). In this case, vehicles of the reference class is clustered into $\Delta n$ sized group. Then, the average spacing $s$ of each vehicle class over the clusters of the reference class is updated at each time step $\Delta t$.\par
For the reference class ($r$) the average spacing is updated following Eq. (11), and the trajectory is updated for according to Eq. (13).

\begin{figure}[!thpb]
  \centering
  \includegraphics[width=3in]{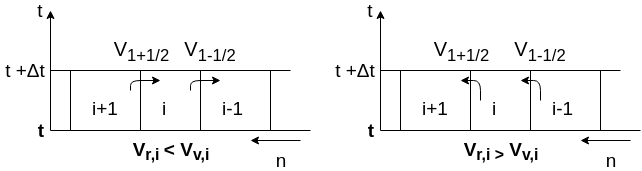}
 \caption{Direction of fluxes through the edges of the cluster}
  \label{fig:fluxdir}
 \end{figure}

The average spacing of the remaining vehicle classes is updated according to:
\begin{equation}
{\left(\frac{s_{r,i}}{s_{c,i}}\right)}^{t+\Delta t}={\left(\frac{s_{r,i}}{s_{c,i}}\right)}^{t} -\frac{\Delta t}{\Delta n}(V_{c,i+1/2}-V_{c,i-1/2}) 
\end{equation}
where $V_{v,i\pm 1/2}$ are the fluxes (speeds) at the cell boundaries. When the speed of the reference class is always higher than the rest classes $(v_r > v_c)$, the direction of the fluxes is to the left. Thus, 
\begin{IEEEeqnarray}{rCl}\label{eq:Ldir_flux}
V_{c,i+1/2}=\frac{v_{r,i}-v_{c,i}}{s_{c,i}}
\IEEEyessubnumber\\
V_{c,i-1/2}=\frac{v_{r,i-1}-v_{c,i-1}}{s_{c,i-1}} 
\IEEEyessubnumber
\end{IEEEeqnarray}
On the other hand, if $(v_r < v_c)$, the direction of the fluxes is to the right (see Fig. \ref{fig:fluxdir}). This suggests that the fluxes should be defined as
\begin{IEEEeqnarray}{rCl} \label{eq:Rdir_flux}
V_{c,i+1/2}=\frac{v_{r,i+1}-v_{c,i+1}}{s_{c,i+1}}
\IEEEyessubnumber\\
V_{c,i-1/2}=\frac{v_{r,i}-v_{c,i}}{s_{c,i}} 
\IEEEyessubnumber
\end{IEEEeqnarray}

However, the flux definition in Eq. (\ref{eq:Rdir_flux}) is restricted to the situation where the fluxes through the edges are non-zero, i.e. $v_{c,i+1}>v_{r,i}$ and $v_{c,i}>v_{r,i-1}$. \par
\begin{figure}[!thpb]
\captionsetup{justification=centering}
\begin{subfigure}{0.45\linewidth}
\includegraphics[width=\textwidth]{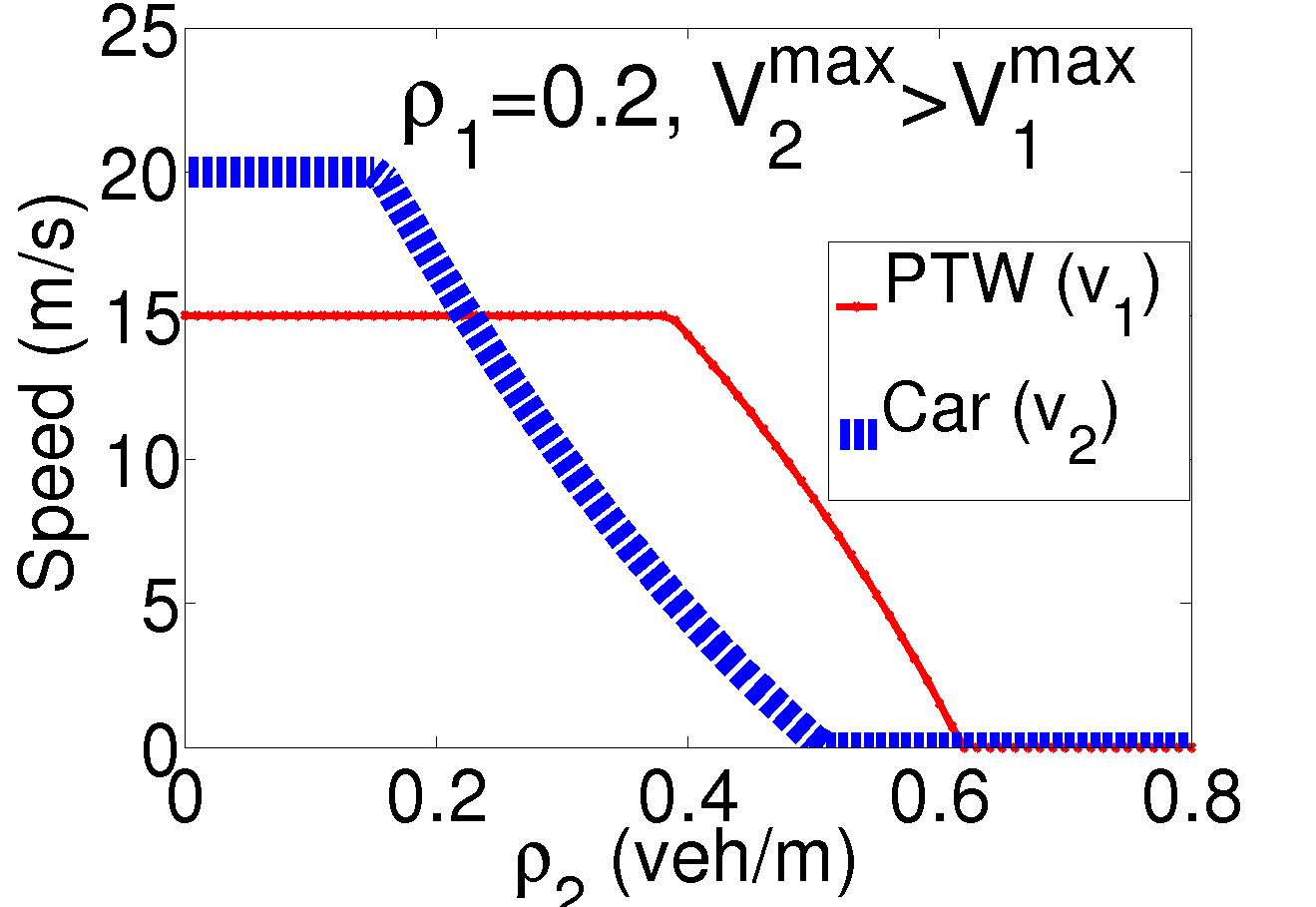}
\caption{}
\label{fig:speeda}
\label{}
\end{subfigure}
\begin{subfigure}{0.45\linewidth}
\includegraphics[width=\textwidth]{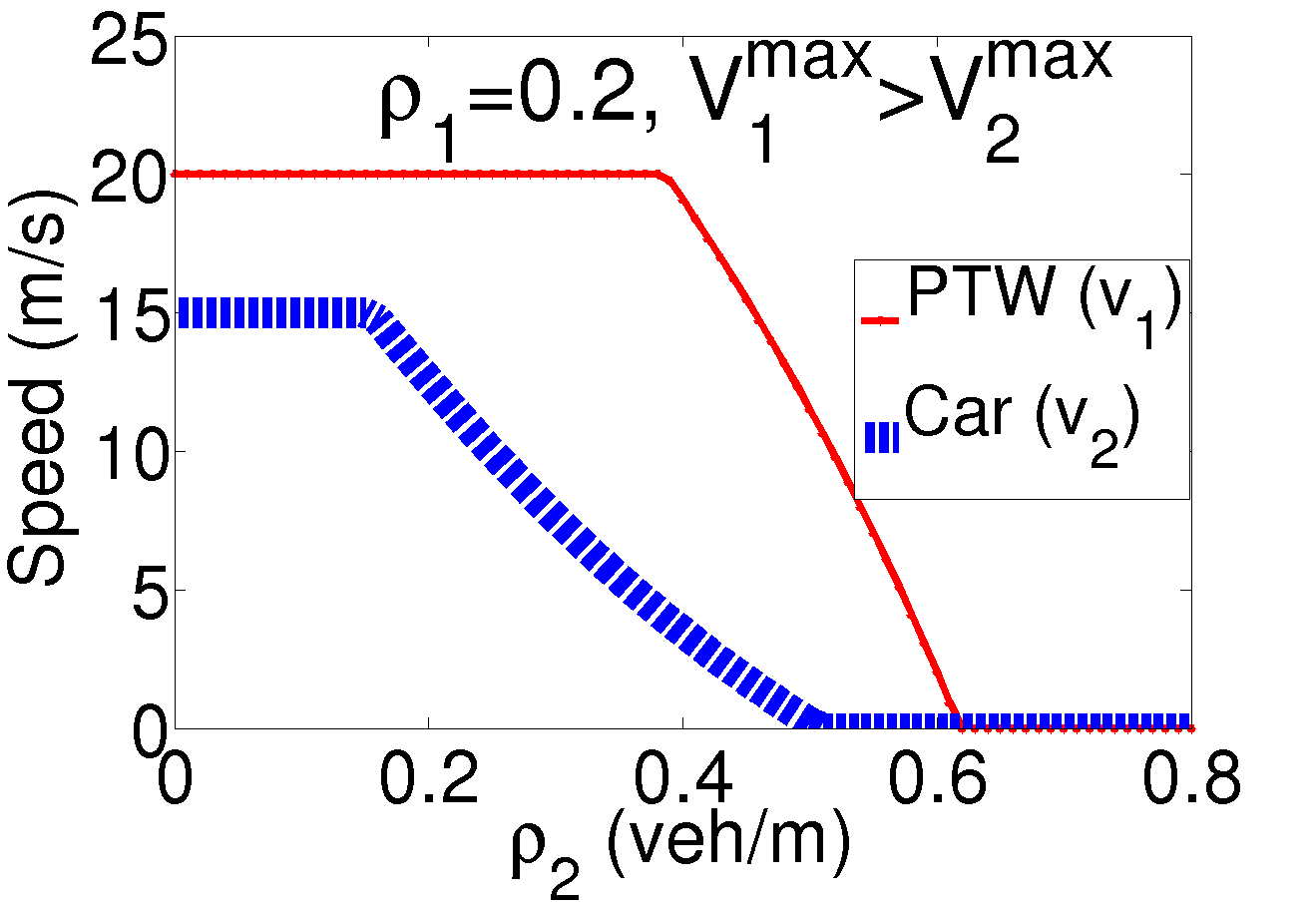}
\caption{}
\label{fig:speedb}
\end{subfigure}
\caption{Speed density relation for (a) free flow speed of Cars is greater PTWs (b) free flow speed of PTWs is greater cars, density of PTWs $\rho_1=0.2 veh/m$}
\end{figure}
For traffic flow consisting of PTWs and cars, if the reference class is PTWs and PTWs have a higher free flow speed than cars (Fig. \ref{fig:speedb}), the flux definition in Eq. (\ref{eq:Ldir_flux}) applies. Nonetheless, if the free flow speed of car is higher than PTWs (Fig. \ref{fig:speeda}), whichever class is the reference class, we have both conditions, $v_r>v_c$ and $v_c>v_s$. For this reason, we give a general definition for the fluxes which applies irrespective of the order of the speeds. \par

If $v_{r,i}^n>v_{c,i}^n$,
\begin{IEEEeqnarray}{lCr}
V_{i+1/2}=\frac{(v_{r,i}-v_{c,i})}{s_{c,i}}
\IEEEyessubnumber\\
V_{i-1/2}=\frac{\max(0,(v_{r,i-1}-v_{c,i-1}))}{(v_{r,i-1}-v_{c,i-1})}\frac{(v_{r,i-1}-v_{c,i-1})}{s_{c,i-1}}\qquad
\IEEEyessubnumber
\end{IEEEeqnarray}

If $v_{r,i}^n<v_{c,i}^n$,
\begin{IEEEeqnarray}{lCr}
V_{i+1/2}=\frac{\max(0,(v_{c,i+1}-v_{r,i}))}{(v_{c,i+1}-v_r^{i}))} \frac{(v_{r,i}-v_{c,i+1})}{s_{c,i+1}}
\IEEEyessubnumber\\
V_{i-1/2}=\frac{\max(0,(v_{c,i}-v_{r,i-1}))}{(v_{c,i}-v_{r,i-1}))} \frac{(v_{r,i-1}-v_{c,i})}{s_{c,i}}
\IEEEyessubnumber
\end{IEEEeqnarray}

\subsection{Follow-the-leader type model from Lagrangian representation}
In continuum flow model, $\Delta n$ can take any positive value. A follow-the-leader type flow is observed when $\Delta n=1$. On the discretization scheme, grouping is done per vehicle classes base, which perfectly works for traffic flows obeying lane discipline. However, when we have two-wheelers which do not respect such an ordered flow, a special treatment is required. The reason is that, in the discretization, clusters of the same vehicle class are not allowed to overlap or occupy the same position. Consequently, the parallel movement of two-wheelers cannot be modeled properly. Thus, we integrate the abreast movement of two-wheelers by introducing sub-lanes. Accordingly, two-wheelers in a sub-lane adhere to the follow-the-leader principle. \par

\begin{figure}[thpb]
  \centering
\includegraphics[width=0.4\textwidth]{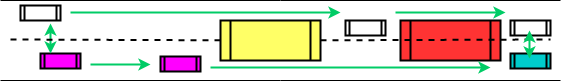}
 \caption{View of vehicle in Lagrangian framework when sublane introduced}
  \label{}
 \end{figure}
\begin{equation}\label{eq:locupdate}
X(i,t+\Delta t)= X(i,t) + \Delta t*V
\end{equation}
The location of the vehicles is updated following Eq. (\ref{eq:locupdate}). Since the macroscopic speed is defined as a function of the free space between vehicles (refer \cite{c8}), the lateral and longitudinal interaction between vehicle classes can be captured. For example, the speed of a PTW depends on the number of vehicles (cars and PTWs on the other sub-lanes) within the space between the leader the follower PTWs and the longitudinal spacing. likewise for cars. With this approach, moving behavior of each vehicle class can be analyzed at a fine-grained level. Further, additional vehicle (or vehicle class) specific rules also can be incorporated, making it a suitable and efficient solution for dealing with cooperative intelligent transport system (C-ITS). 
\FloatBarrier
\section{Numerical results and discussion}
To test the validity and accuracy of the proposed discretization scheme, we compare the numerical results obtained with Eulerian approach and the two Lagrangian methods. For the simulation experiment, the parameters in Table \ref{table:setting} are used.
\begin{table}[h]
\caption{simulation settings}
\label{table:setting}
\begin{center}
\begin{tabular}{l c}
 \hline \\
 Maximum speed of cars & 15m/s\\
 Maximum speed of PTWs & 20m/s\\
 vehicle cluster size & 7.5 vehicles\\
 Time step & 0.125 s\\
 Space steps (Eulerian) & 10m\\
 Road length & 3000m\\
 lane width & 3.5 m\\
 Number of lanes & 1\\
 Simulation time & 45s\\
 \hline
 \end{tabular}
 \end{center}
\end{table}
Lax-friedrich discretization scheme is employed to solve the Eulerian conservation equations. We assume identical initial densities for the two vehicle classes, cars ($\rho_2$) and PTWs ($\rho_1$), where $\rho_1=\rho_2=0.15 veh/m$, for $x \in [0, 1400m]$ and $\rho_1=\rho_2=0.3 veh/m$, otherwise. \par
\begin{figure}[!htbp]
\centering
\begin{subfigure}{\linewidth}
\includegraphics[width=\textwidth]{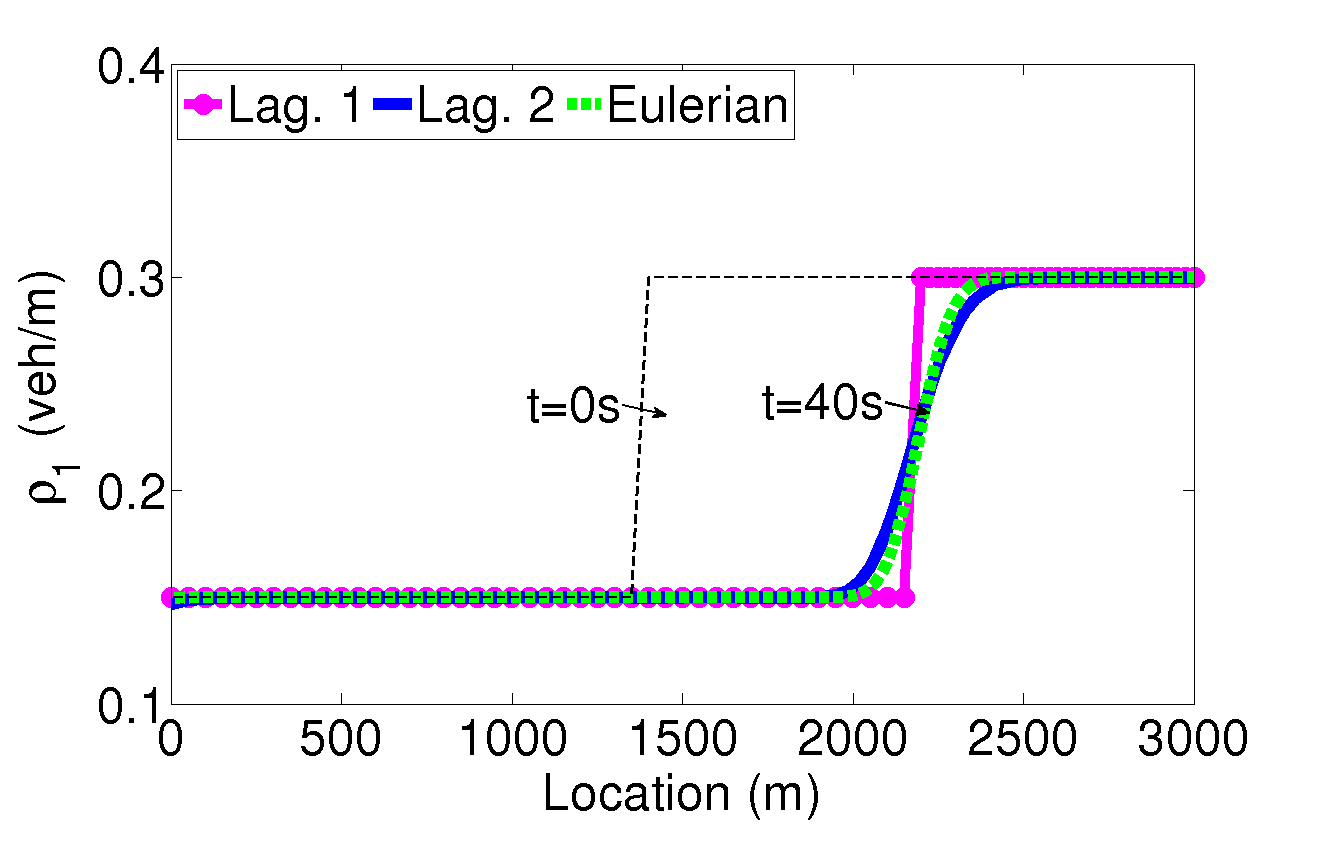}
\caption{PTW density wave}
\label{fig:lageul_ptw}
\end{subfigure}
\begin{subfigure}{\linewidth}
\includegraphics[width=\textwidth]{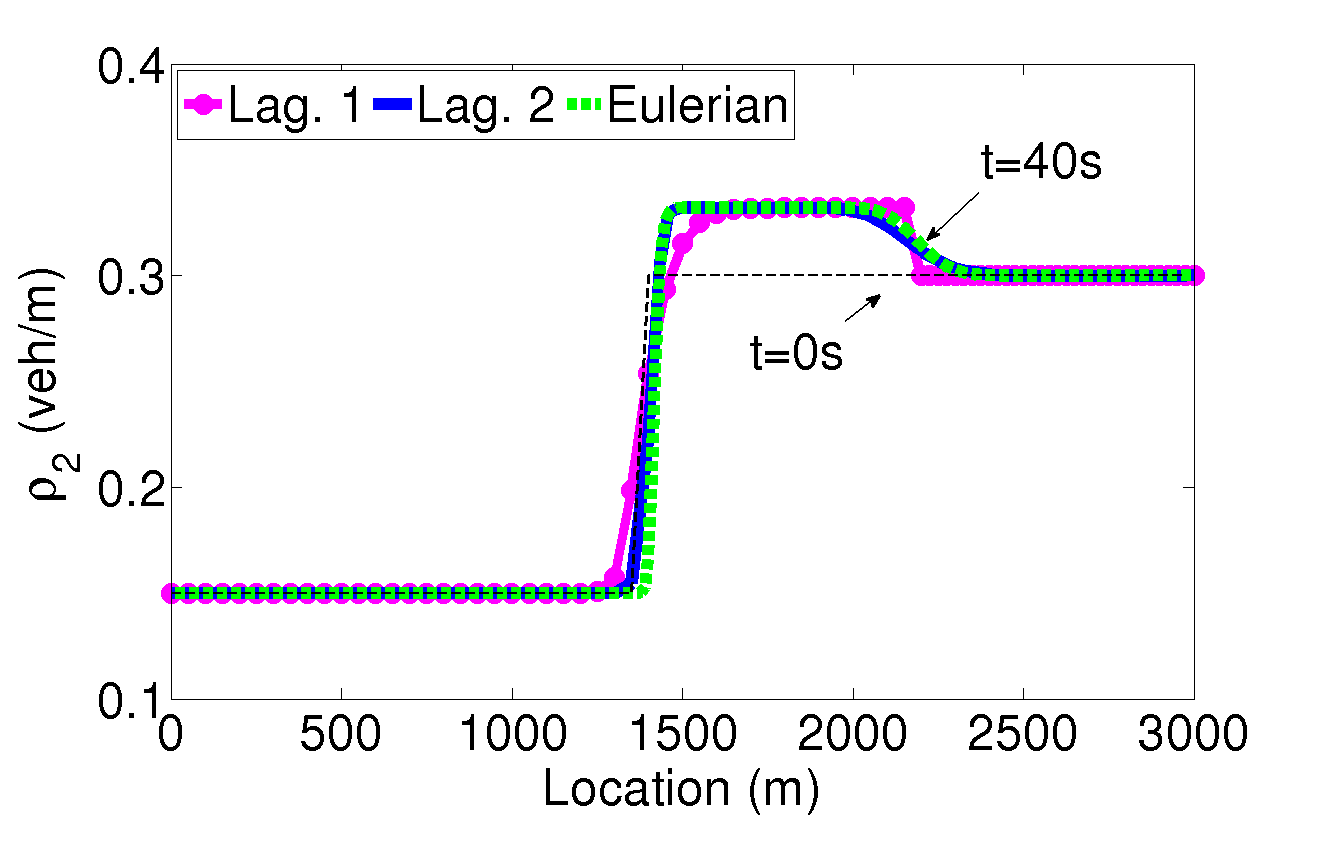}
\caption{Car density wave}
\label{fig:lageul_car}
\end{subfigure}
\caption{Lagrangian when PTW is the reference class ( Lag. 1) and cars is the reference class (Lag. 2) vs Eulerian}
\label{fig:lageul}
\end{figure}

The evolution of the initial density as described by Eulerian and Lagrangian approach is presented in Fig. \ref{fig:lageul}. For the Lagrangian approach we considered two cases by changing the reference class. Lag. 1 stands for the result when PTWs are the reference class and Lag. 2 stands for the results when cars are the reference class. In this case, the fundamental diagram takes the shape in \ref{fig:speedb}.\par


\FloatBarrier
The density wave of PTWs and cars at time $t=40s$ are depicted in Fig. \ref{fig:lageul_car} and \ref{fig:lageul_ptw}, respectively. As can be seen, the results are close to each other except the difference at the upstream and downstream shock fronts. With this, we can prove the validity of the proposed discretization scheme for the case where the slower vehicle is the reference vehicle class (shown by Lag. 2).\par

\begin{figure}[!htbp]
\centering
\begin{subfigure}{\linewidth}
\includegraphics[width=\textwidth]{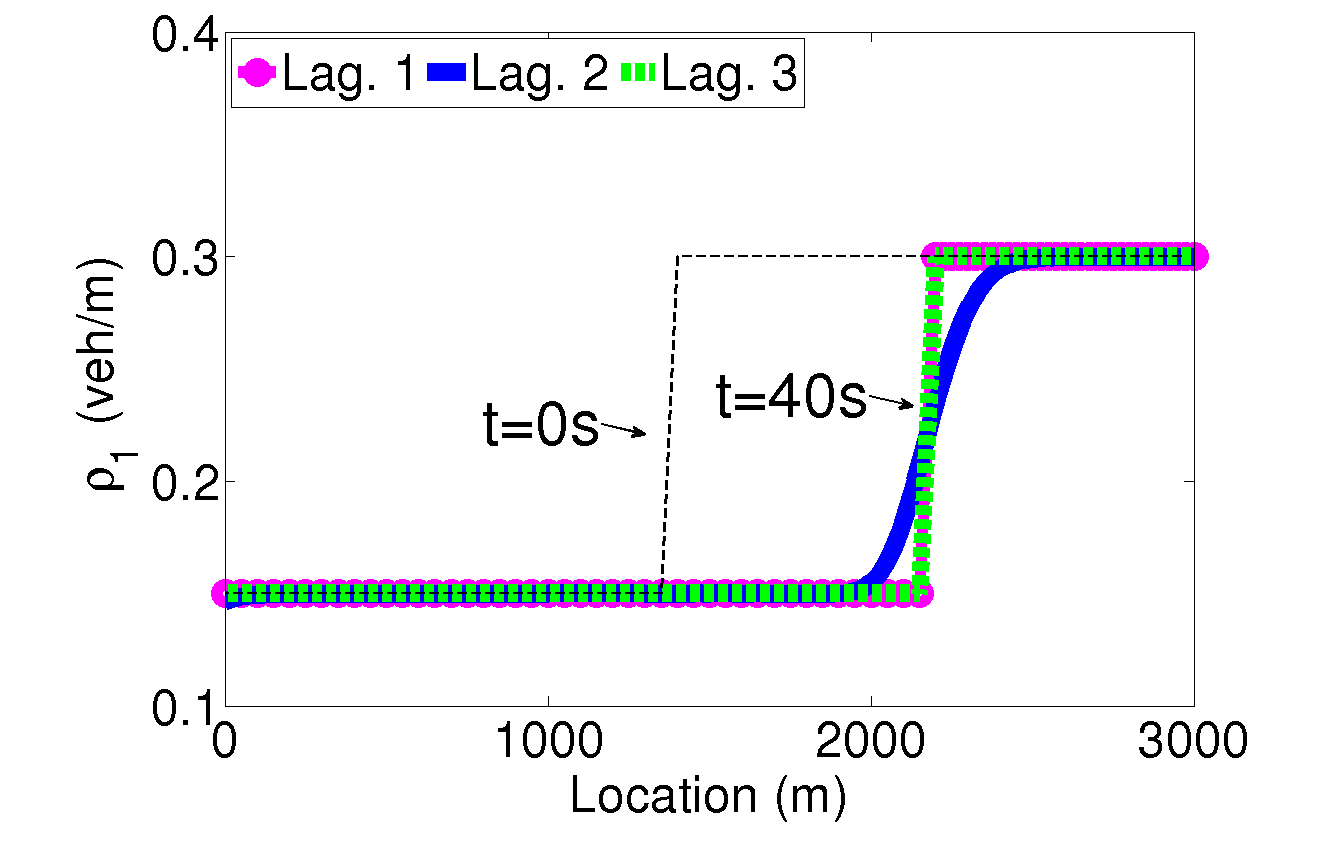}
\caption{PTW density wave}
\label{fig:laglag_ptw}
\end{subfigure}
\begin{subfigure}{\linewidth}
\includegraphics[width=\textwidth]{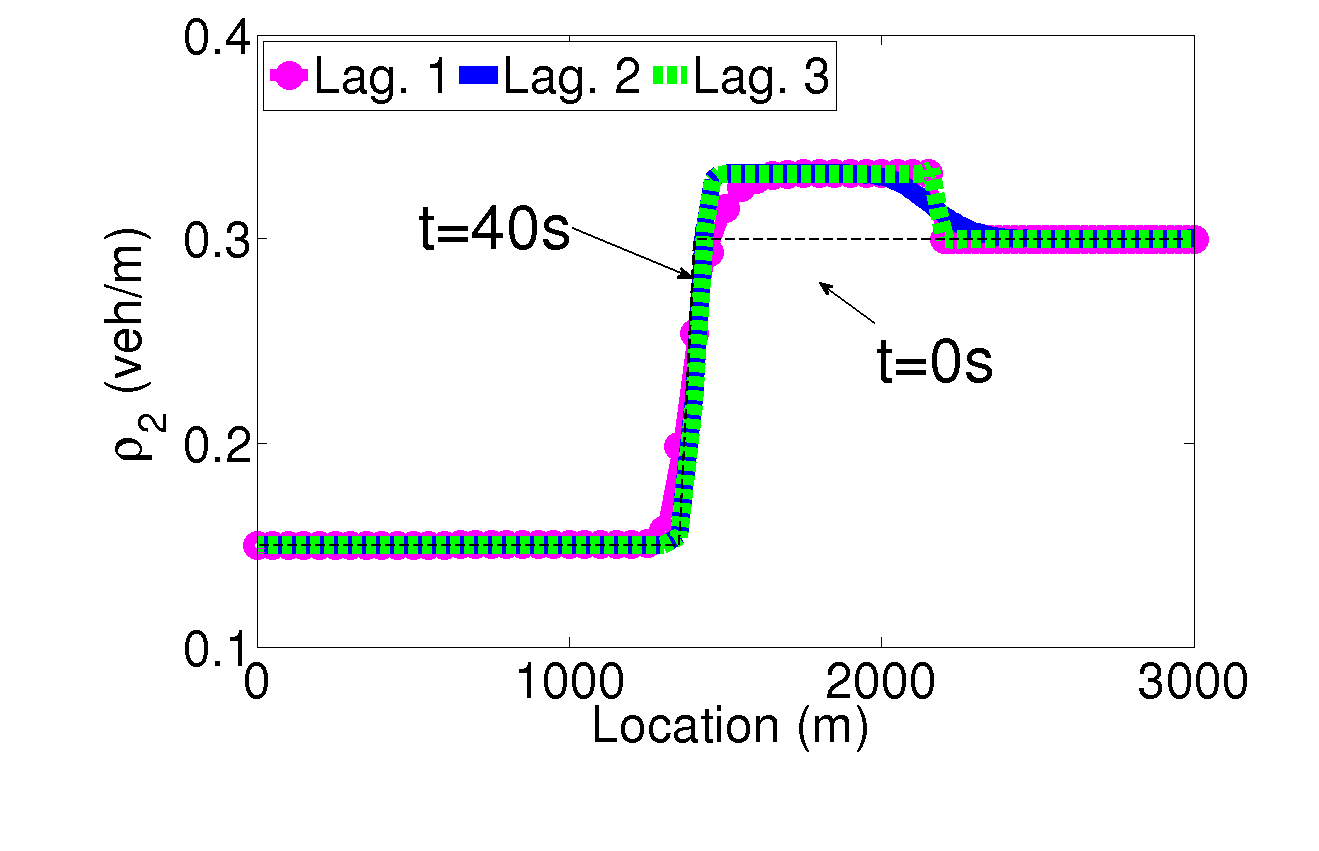}
\caption{Car density wave}
\label{fig:laglag_car}
\end{subfigure}
\caption{Comparison between the two Lagrangian methods, (Lag. 1, Lag.2 ) for method 2, (Lag. 3) for method 1.}
\label{fig:laglag}
\end{figure}

Furthermore, the comparison of the two Lagrangian methods is presented in Fig. \ref{fig:laglag}. The density waves for cars and PTWs shown in Figs. \ref{fig:laglag_car} and \ref{fig:laglag_ptw} illustrate that method 1 (Lag. 3) produces a more accurate result than method 2 (Lag. 1 and Lag. 2) of the Lagrangian approach. Specifically, at the  high density to low density and low density to high density transition points numerical error are observed for the case of Lag. 1 and Lag. 2. \par
\begin{figure}[!htbp]
\centering
\includegraphics[width=3.5in]{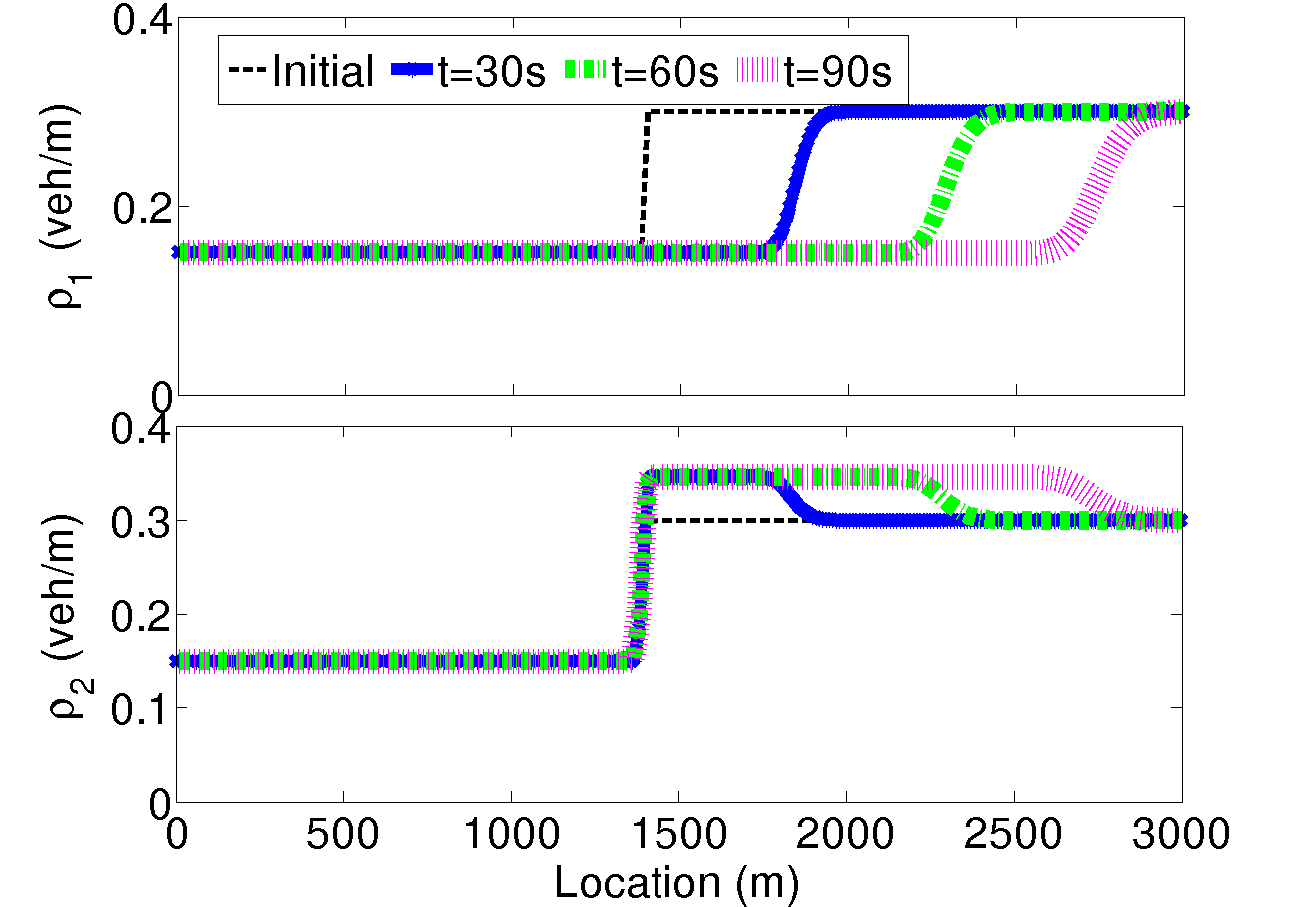}
\caption{Density waves of PTWs (upper subplot) and cars (lower subplot), the numerical scheme defined for the general case of method 2 is applied}
\label{fig:gen_flux}
\end{figure}

We also test the proposed numerical scheme, i.e. the definition of the fluxes at the boundary, for method 2 Lagrangian representation. For this experiment, we consider the fundamental digram in Fig. \ref{fig:speeda}, and the maximum speed of cars$=20 m/s$ and the maximum speed of PTWs$=15 m/s$. The rest simulation parameters and the initial density are identical to the the previous experiments. The evolution of cars and PTWs density waves is shown in Fig. \ref{fig:gen_flux}. According to the result obtained, the evolution is correctly described by the proposed scheme.\par

\begin{figure}[!htbp]
\centering
\includegraphics[width=3.5in]{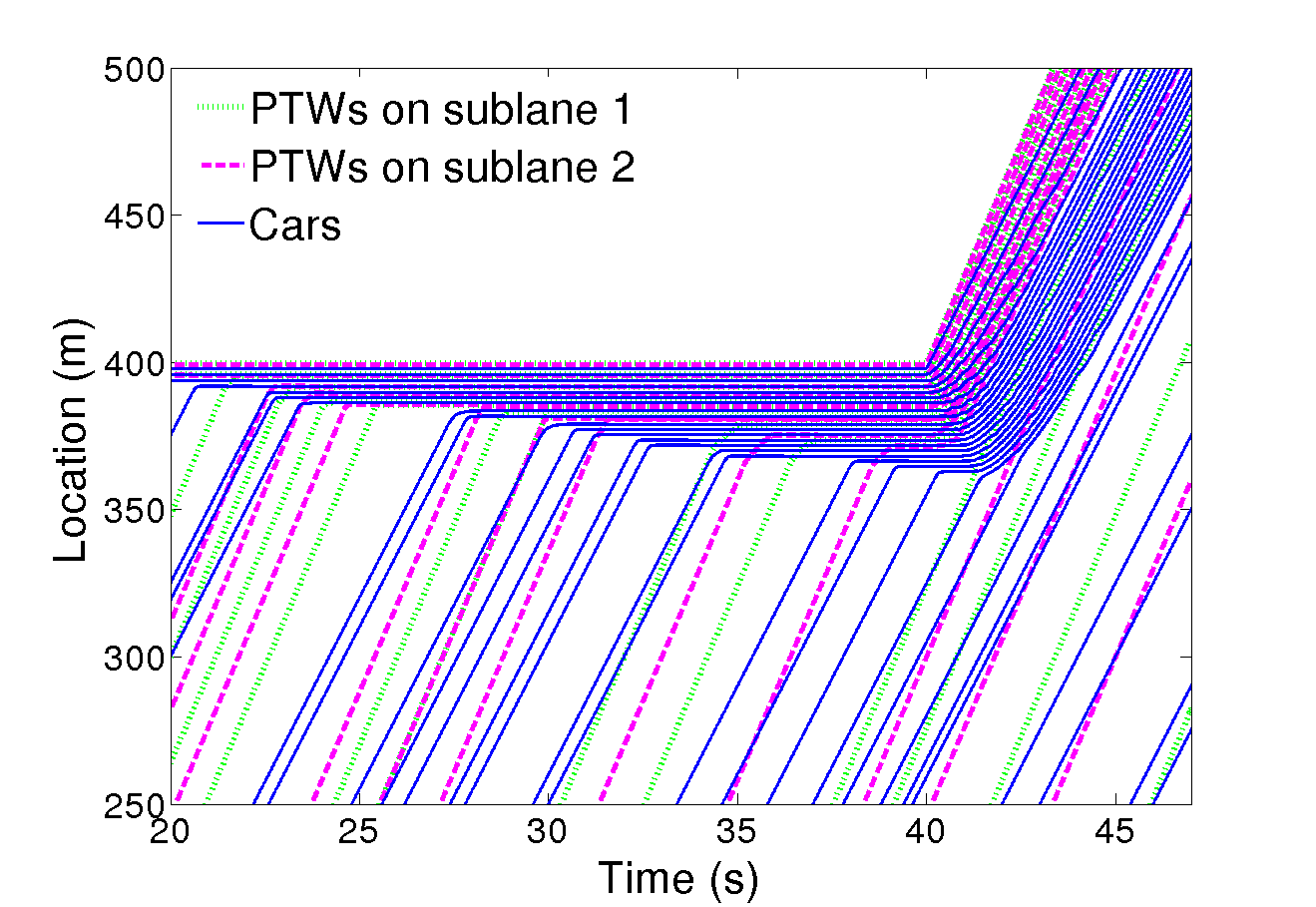}
\caption{Trajectories of vehicles, two sub-lanes for PTWs}
\label{fig:tra}
\end{figure}

For the case $\Delta n=1$ the trajectory of the vehicles on the space-time plane is presented in Fig. \ref{fig:tra}. To track the interaction between vehicle classes at different traffic situations, a traffic light is located at $400m$ which stays red for the period $t \in [0,40s]$. PTWs have two sub-groups (sub-lanes) and the clustering of each sub-groups is done separately. As can be observed from overlapping trajectories of PTWs, by introducing sub-lanes the side by side movement of PTWs, in the same lane, can be reproduced (see the trajectory of vehicles departing the queue).  
\section{CONCLUSIONS}
Lagrangian formulation gives accurate representation, permits to study various traffic feature and is applicable for the current traffic state estimation schemes. Due to these benefits, Lagrangian representation is preferred over the Eulerian one. In this paper, we formulate the multiclass LWR model for a traffic flow consisting of PTWs and cars in Lagrangian coordinates. \par
 We proposed a numerical scheme taking into account the peculiar features observed in mixed flow of cars and PTWs. The validity of proposed method checked through simulation experiments. Accordingly to the result, our numerical scheme can produce a valid results. Moreover, the simulation results shows that the Lagrangian representation outperform over the Eulerian representation in terms of accuracy. The possibility of tracking the trajectory of each vehicles in Lagrangian representation, facilitate the investigation of different traffic phenomena for C-ITS applications.
 



\FloatBarrier
\section*{ACKNOWLEDGMENT}
This work was funded by the French Government (National Research Agency, ANR) through the “Investments for the Future” Program reference \#ANR-11-LABX-0031-01.
EURECOM acknowledges the support of its industrial members, namely, Orange, BMW Group, SAP, Monaco Telecom, Symantec, IABG.



\begin{thebibliography}{99}
\bibitem{c1} Leclercq, L., Laval, J. A. (2009). A multiclass car-following rule based on the LWR model. In Traffic and Granular Flow’07 (pp. 151-160). Springer, Berlin, Heidelberg.
\bibitem{c2} Leclercq, L., Laval, J. A., Chevallier, E. (2007). The Lagrangian coordinates and what it means for first order traffic flow models. In Transportation and Traffic Theory 2007. Papers Selected for Presentation at ISTTT17.
\bibitem{c3} Newell, G. F. (2002). A simplified car-following theory: a lower order model. Transportation Research Part B: Methodological, 36(3), 195-205.
\bibitem{c4}van Wageningen-Kessels, F., van Lint, H., Hoogendoorn, S., Vuik, K. (2010). Lagrangian formulation of multiclass kinematic wave model. Transportation Research Record: Journal of the Transportation Research Board, (2188), 29-36.
\bibitem{c5}Newell, G. F. (1993). A simplified theory of kinematic waves in highway traffic, part I: General theory. Transportation Research Part B: Methodological, 27(4), 281-287.
\bibitem{c6}Lighthill, M. J., Whitham, G. B. (1955, May). On kinematic waves. I. Flood movement in long rivers. In Proceedings of the Royal Society of London A: Mathematical, Physical and Engineering Sciences (Vol. 229, No. 1178, pp. 281-316). The Royal Society.
\bibitem{c7}Richards, P. I. (1956). Shock waves on the highway. Operations research, 4(1), 42-51.
\bibitem{c8}Gashaw, S.M., Goatin, P. and H\"{a}rri, J., (2018, April). Modeling and analysis of mixed flow of cars and powered two wheeelers, Transportation Research Part C: Emerging Technologies, Volume 89, Pages 148-167, ISSN 0968-090X,
\bibitem{c9} Yuan, Y., F. L. M. van Wageningen-Kessels, J. W. C. van Lint, S. P. Hoogendoorn
(2011) Two modeling and discretization choices for Lagrangian multi-class first-
order traffic flow model and their related (dis-)advantages, in: The Ninth Interna-
tional Conference on Traffic and Granular Flow, 2011.
\bibitem{c10} Ni, D. (2007). Determining traffic-flow characteristics by definition for application in ITS. IEEE Transactions on Intelligent Transportation Systems, 8(2), 181-187.
\bibitem{c11} Yuan, Y., Van Lint, J. W. C., Wilson, R. E., van Wageningen-Kessels, F., Hoogendoorn, S. P. (2012). Real-time Lagrangian traffic state estimator for freeways. IEEE Transactions on Intelligent Transportation Systems, 13(1), 59-70.
\bibitem{c12} Seo, T.,  Kusakabe, T. (2015). Probe vehicle-based traffic state estimation method with spacing information and conservation law. Transportation Research Part C: Emerging Technologies, 59, 391-403.
\bibitem{c13}Leclercq, L. (2007). Hybrid approaches to the solutions of the “Lighthill–Whitham–Richards” model. Transportation Research Part B: Methodological, 41(7), 701-709.
\bibitem{c15} Yuan, K., Knoop, V. L., Hoogendoorn, S. P. (2017). A kinematic wave model in Lagrangian coordinates incorporating capacity drop: Application to homogeneous road stretches and discontinuities. Physica A: Statistical Mechanics and its Applications, 465, 472-485.
\end{thebibliography}
\end{document}